\documentclass[lettersize,journal]{IEEEtran}
\usepackage{amsmath}
\usepackage{graphicx}
\usepackage{amsmath}
\usepackage{epstopdf} 
\usepackage{mathptmx}      
\usepackage[utf8]{inputenc}
\usepackage{bm}
\usepackage{amsfonts}
\usepackage{graphics}
\usepackage{amssymb}
\usepackage{algorithmicx}
\usepackage{algorithm}
\usepackage{algpseudocode}
\usepackage{multirow}
\usepackage{xcolor}
\usepackage{url}
\usepackage[stable]{footmisc}
\usepackage{amssymb}
\begin{document}
	
	\title{An iterative Jacobi-like algorithm to compute a few sparse approximate eigenvectors}
	
	\author{Cristian Rusu
		\thanks{C. Rusu is with the Faculty of Mathematics and Computer Science, Bucharest, 030018, Romania. E-mail: cristian.rusu@unibuc.ro}}

\markboth{Journal of \LaTeX\ Class Files,~Vol.~14, No.~8, August~2021}%
{Shell \MakeLowercase{\textit{et al.}}: A Sample Article Using IEEEtran.cls for IEEE Journals}


\maketitle

\begin{abstract}
	In this paper, we describe a new algorithm that approximates the extreme eigenvalue/eigenvector pairs of a symmetric matrix. The proposed algorithm can be viewed as an extension of the Jacobi eigenvalue method for symmetric matrices diagonalization to the case where we want to approximate just a few extreme eigenvalues/eigenvectors. The method is also particularly well-suited for the computation of sparse approximations of the eigenvectors. In fact, we show that in general, our method provides a trade-off between the sparsity of the computed approximate eigenspaces and their accuracy. We provide theoretical results that show the linear convergence of the proposed method. Finally, we show experimental numerical results for sparse low-rank approximations of random symmetric matrices and show applications to graph Fourier transforms, and the sparse principal component analysis in image classification experiments. These applications are chosen because, in these cases, there is no need to perform the eigenvalue decomposition to high precision to achieve good numerical results.
\end{abstract}


\section{Introduction}
Low-rank eigenvalue decompositions (EVD) of matrices~\cite{DecompositionsEig} are one of the most important factorizations in numerical linear algebra with many applications in applied mathematics, machine learning, and signal processing.

For full symmetric matrices, one of the most popular ways to build eigendecompositions is the QR algorithm~\cite{10.1093/comjnl/4.4.332,814656}. While this algorithm efficiently produces the complete eigenfactorization, it can be numerically expensive to run if the input matrix is very large or has some special structure that is destroyed by the algorithm. For sparse and/or large scale matrices, several methods have been proposed in the literature which compute just a few eigenvalues/eigenvectors pairs of interest (either a few of the extreme ones or a few near a given fixed eigenvalue $\lambda$): iterative subspace methods \cite[Chapter~7.3]{Golub1996} with Rayleigh-Ritz acceleration~\cite{10.1007/BF02219773}, Lanczos methods~\cite{SIMON1984101} with restarts~\cite{10.1007/BFb0018526}, Jacobi-Davidson methods~\cite{10.2307/2653109,Ravibabu2020Dec}, Rayleigh quotient~\cite{LONGSINE1980195} and trace minimization~\cite{SAMEH2000155}.


Introduced in 1848, the Jacobi method for the diagonalization of a symmetric matrix~\cite{JacobiProcess} is a conceptually simple yet effective method to find all the eigenvalues/eigenvectors of a symmetric matrix. {\color{black}The method iteratively transforms the original matrix by setting to zero the largest absolute value off-diagonal element in the matrix until we are left only with the diagonal elements, i.e., the eigenvalues (up to some precision)}.

The numerical computation of eigenvectors is still of great interest as there is much research underway to perform faster~\cite{9113480}, structured~\cite{8962381}, and more robust decompositions~\cite{9406364} with modern applications~\cite{LIU2021107931}. In recent years, driven mainly by applications in machine learning and signal processing, researchers have started to explore algorithms and theory for the construction of {\color{black}sparse approximation of eigenvectors}. In the era of big data, it is convenient to find eigenspaces that are easier to store (due to sparsity) and interpretable. To this end, the sparse principal component analysis (SPCA) method was introduced in~\cite{doi:10.1198/106186006X113430} to build {\color{black}sparse eigenvectors}. As SPCA involves a sparsity constraint, it is NP-hard to solve in general~\cite{6658871} and, therefore, there is continuous research to find better algorithms and tighter guarantees to find sparse principal components. There are several strategies proposed to solve the SPCA problem: a generalization of the power method \cite{10.5555/1756006.1756021} that allows $\ell_1$ sparsity regularization or a truncated variant of the power method \cite{10.5555/2567709.2502610}, the work in~\cite{doi:10.1198/1061860032148,doi:10.1198/106186006X113430,SMALLMAN2018443} added an $\ell_1$ regularization to achieve sparsity, a direct formulation based on semidefinite programming~\cite{10.2307/20453990}, greedy methods for SPCA~\cite{NIPS2005_d4784467} and iterative thresholding methods applied to the SVD~\cite{SHEN20081015,10.1093/biostatistics/kxp008}. More recent work uses branch-and-bound strategies~\cite{Berk2019}, bipartite matching~\cite{NIPS2015_2b8a6159}, low-rank updates~\cite{papailiopoulos2013sparse} and coordinate descent~\cite{Beck2016}.

In this paper, we propose an extension of the Jacobi method that computes only a few eigenvectors and is therefore well suited to build $p$-rank approximations of symmetric matrices. The method is also particularly well suited when calculating {\color{black}sparse eigenvectors}. We note that we are looking for approximate decompositions and therefore we will highlight a trade-off between the sparsity or complexity of the eigenspaces and the accuracy of the results we compute. We provide theoretical results that show linear convergence of the proposed method and give insights into the sparsity levels reached. We highlight the fact that while some computed eigenspaces might not be sparse, the proposed factorization allows for their numerically efficient manipulations (projection operations).

The paper is organized as follows: Section 2 gives an outline of the Jacobi method for matrix diagonalization, Section 3 describes the main ideas of the paper and the proposed algorithm together with a convergence analysis, and finally, Section 4 gives experimental numerical results.

{\color{black}\noindent\textbf{Notation.} We denote $\mathbf{x} \in \mathbb{R}^n_{+}$ a vector of size $n$ whose entries are nonnegative and similarly $\mathbf{x} \in \mathbb{R}^n_{-}$ is vector of size $n$ whose entries are non-positive, the set $\mathbb{N}$ is the set of natural numbers, $\mathbf{I}_n$ is the identity matrix of size $n \times n$, $\mathbf{1}_{n \times m}$ is the matrix of size $n \times m$ with all elements one. Also, $a \leftarrow b$ reads as ``$a$ takes the value $b$'' while the function diag($\mathbf{x}$) takes as input a vector $\mathbf{x}$ of size $n$ and returns a matrix $\mathbf{X}$ of size $n \times n$ whose diagonal are the elements of the vector $\mathbf{x}$.}

\section{The Jacobi method for the diagonalization of a symmetric matrix}
The Jacobi eigenvalue method is an iterative procedure that computes all the eigenvalues/eigenvectors of a real symmetric matrix. Starting from the given symmetric matrix $\mathbf{S}^{(0)}$ {\color{black}the method computes} a series of updates like $\mathbf{S}^{(q)} \leftarrow \mathbf{J}_{i_q j_q}^T\mathbf{S}^{(q-1)} \mathbf{J}_{i_q j_q}$ where $\mathbf{J}_{i_qj_q}$ is a Jacobi rotation such that the entry $(i_q,j_q)$ of the matrix $\mathbf{S}^{(q)}$ is zeroed. At step $q$, we choose to zero the largest absolute value entry in the matrix, i.e., $(i_q,j_q) = \underset{v > t}{\arg \max}\ |S_{tv}^{(q)}|$. This guarantees that the iterates converge to a diagonally dominant matrix, i.e., $\lim_{q \to \infty} \mathbf{S}^{(q)} = \text{diag}(\mathbf{\lambda})$ ({\color{black}$\mathbf{\lambda} \in \mathbb{R}^{n}$ is a vector and $\bm{\Lambda} = \text{diag}(\mathbf{\lambda})$ is the diagonal of eigenvalues}). The method converges linearly at first and quadratically after a certain number of updates~\cite{10.2307/2098731}. {\color{black}Moreover, it can be implemented efficiently on multi-core and multi-processor computing architectures ~\cite{BECKA2002243,doi:10.1142/S0129626415500036,Shroff1990APA}.}

\section{The proposed method}

In this section, we modify the Jacobi method to approximately compute only a few eigenvalues/eigenvectors. Similar to the Jacobi method, we are concerned with two questions:
\color{black}
\begin{enumerate}
	\item[(A)] how to choose the pivot indices $(i
	_q,j_q)$?
	\item[(B)] {\color{black}what similarity transformation to apply on the current matrix $\mathbf{S}^{(q)}$?}
\end{enumerate}

{\color{black}Our solution is based on the following three steps: 1) we write the eigenvalue problem as a least-squares problem whose solution is a $p$-rank eigenspace; 2) to solve this least-square problem, we introduce basic building blocks for which we have closed-form solutions; and 3) we show that the eigenspace can be written as a product of these basic building blocks. Then, we give an efficient algorithm to construct an approximate eigenspace and provide theoretical convergence results that express the quality of the result.}

\subsection{The main building blocks}

In the following proposition, we provide the main result for the basic building blocks on which we based our algorithm.

{\color{black}\noindent \textbf{Proposition 1.} Let $\mathbf{\alpha} \in \mathbb{R}^p, 1 \leq p \leq n$ and let $\mathbf{S} \in \mathbb{R}^{n \times n}$ be a symmetric matrix of order $n$. Let $\mathbf{S} = \mathbf{U} \text{diag}(\mathbf{\lambda}) \mathbf{U}^T$ be its eigenvalue decomposition with the vector $\mathbf{\lambda} \in \mathbb{R}^n$ and $\lambda_1 \geq \lambda_2 \geq \dots \geq \lambda_n$. Then the following optimization problem in the variable $\mathbf{\overline{U}} \in \mathbb{R}^{n \times n}$}:
\begin{equation}
	\underset{\mathbf{\overline{U}},\ \mathbf{\overline{U}}^T \mathbf{\overline{U}} = \mathbf{I}_n }{\text{minimize}} \ \| \begin{bmatrix} \text{diag}(\mathbf{\alpha}) & \mathbf{0}_{p \times (n-p)} \\ \mathbf{0}_{(n-p) \times p} & \mathbf{0}_{(n-p) \times (n-p)} \end{bmatrix} - \mathbf{\overline{U}}^T\mathbf{S}\mathbf{\overline{U}} \|_F^2,
	\label{eq:optimization}
\end{equation}
has the following solutions and objective function values:
\begin{itemize}
	\item $\mathbf{\alpha} \in \mathbb{R}_+^{p}$ is sorted increasingly: $\mathbf{\overline{U}}$ contains the $p$ eigenvectors associated with the largest $p$ eigenvalues of $\mathbf{S}$ and the objective function value is given by $\sum_{t=1}^p (\alpha_t-\lambda_t)^2 + \sum_{t=p+1}^n \lambda_t^2$;
	\item $\mathbf{\alpha} \in \mathbb{R}_-^{p}$ is sorted decreasingly: $\mathbf{\overline{U}}$ contains the $p$ eigenvectors associated with the smallest $p$ eigenvalues of $\mathbf{S}$ and the objective function value is given by $\sum_{t=1}^{n-p} \lambda_t^2 + \sum_{t=1}^p (\alpha_{p-t+1} - \lambda_{n-p+t})^2$.
\end{itemize}
In both cases, the columns of $\mathbf{\overline{U}}$ are stored in the decreasing order of their associated eigenvalues.$\hfill \blacksquare$

{\color{black}Herein, we denote with $\mathbf{A}$ a diagonal matrix of size $n \times n$ such that $A_{ii} = \alpha_i$ for $1 \leq i \leq p$ and $A_{ii} = 0$ otherwise. While Proposition 1 applies for general $\mathbf{\overline{U}}$, in the following we use the result for a special structured case where only two coordinates (out of the $n$) are non-zero, i.e., essentially we apply Proposition 1 on two dimensions. Therefore, we now consider the particular case when 
	\begin{equation}
		\mathbf{\overline{U}} =
		\mathbf{G}_{ij} = \begin{bmatrix} \mathbf{I}_{i-1} & & &\\
			& * & & * \\
			& & \mathbf{I}_{j-i-1} & \\
			& * & & * \\
			& & & & \mathbf{I}_{n-j} \\
		\end{bmatrix},
		\label{eq:theG}
\end{equation}}
{\color{black}with $\mathbf{\widetilde{G}}_{ij} = \begin{bmatrix} * & * \\ * & * \end{bmatrix} \in \mathbb{R}^{2 \times 2}$ 
	such that $\mathbf{\widetilde{G}}_{ij}^T \mathbf{\widetilde{G}}_{ij} = \mathbf{I}_2$. The matrix $\mathbf{G}_{ij}$ differs from $\mathbf{I}_n$ only in positions $(i,i), (i,j), (j,i)$ and $(j,j)$ where $\mathbf{\widetilde{G}}_{ij}$ is placed.} 

{\color{black}\noindent \textbf{Proposition 2.} Assume that $\mathbf{\overline{U}}$ in~\eqref{eq:optimization} has the structure from \eqref{eq:theG}. Then, the objective function of~\eqref{eq:optimization} is minimized when the non-trial part of $\mathbf{G}_{ij}$ is set to $\mathbf{V}$ from the $2 \times 2$ eigenvalue decomposition of}
\begin{equation}
	\mathbf{\widetilde{S}}_{ij} = \begin{bmatrix}
		S_{ii} & S_{ij} \\ S_{ji} & S_{jj}
	\end{bmatrix} = \mathbf{VDV}^T.
	\label{eq:theGij}
\end{equation}
{\color{black}In this case, the minimum of~\eqref{eq:optimization} is:
	\begin{equation}
		\min \| \mathbf{A} - \mathbf{G}_{ij}^T \mathbf{S} \mathbf{G}_{ij} \|_F^2  =  \|\mathbf{\alpha} \|_2^2 + \| \mathbf{\lambda} \|_2^2 - 2\sum_{t=1}^p  \alpha_t S_{tt} - 2 \mathcal{C}_{ij},
		\label{eq:theoptimizationproblem}
	\end{equation}
	\begin{equation}
		\mathcal{C}_{ij} \! =\! \!
		\begin{cases}
			\! \! \frac{1}{2} (A_{jj} - A_{ii})(S_{ii} - S_{jj}  +\! \! \! \sqrt{(S_{ii} - S_{jj})^2  + 4S_{ij}^2} ), \text{if } A_{ii} \leq A_{jj},\\
			\! \! \frac{1}{2} (A_{jj} - A_{ii})(S_{ii} - S_{jj}  -\! \! \!  \sqrt{(S_{ii} - S_{jj})^2  + 4S_{ij}^2} ), \text{if } A_{ii} > A_{jj}.
		\end{cases}
		\label{eq:badcij}
	\end{equation}
	for $1 \leq i < j \leq n$ and $\mathcal{C}_{ij} = 0$ for all $1 \leq p < i < j \leq n$.}$\hfill\blacksquare$

{\color{black}We would like to note that $\mathbf{G}_{ij}$ as defined in~\eqref{eq:theGij} is a square orthonormal matrix (orthogonal and columns/rows have unit $\ell_2$ norm) and includes as a special case the widely used Givens rotation. Considering both rotations and reflectors is one of the key ideas of this work as it allows us to use Procrustes formulas to reach the result of Proposition 2.}

{\color{black}It is clear from Proposition 1 that it would be ideal to have the true eigenvalues of $\mathbf{S}$ for which we want to recover the eigenvectors. Unfortunately, these are not available in general. Therefore we propose two ways of choosing $\mathbf{\alpha}$: i) either $\alpha_i \in \{\pm 1\}$ (as per Proposition 1, {\color{black}to recover largest and smallest eigenvalues/eigenvector pairs}, respectively) or ii) $\mathbf{\alpha}$ is a decreasing/increasing series. In the first case, if the entries of $\mathbf{\alpha}$ are equal then $\mathcal{C}_{ij} = 0$ for all $i,j \leq p$ and therefore we will have for large $k$ that $\mathbf{S}^{(k)} \approx \begin{bmatrix}
		\mathbf{S}_1 & \mathbf{0}_{p \times (n-p)} \\
		\mathbf{0}_{(n-p) \times p} & \mathbf{S}_2
	\end{bmatrix}$ where the two blocks $\mathbf{S}_1 \in \mathbb{R}^{p \times p}$ and $\mathbf{S}_2 \in \mathbb{R}^{(n-p) \times (n-p)}$ are symmetric and they split the $n$ eigenvalues of $\mathbf{S}$. {\color{black}Note that $A_{jj} - A_{ii}$ is zero in two situations: either $A_{jj} = A_{ii}$ (which is true when $i, j \leq p$ because it is equivalent to $\alpha_j = \alpha_i$) or when $A_{jj} = A_{ii} = 0$ (always true when $i,j > p$). This means that $\mathcal{C}_{ij}$ can be non-zero only when $i \leq p$ and $j >p$ and therefore the transformations $\mathbf{G}_{ij}$ can only compute a linear combination of the eigenvectors (in some applications, like principal component analysis~\cite{SEGHOUANE2018311,9429896,LIM2021108096}, this is sufficient). These scores $\mathcal{C}_{ij}$, which are non-zero only when $i \leq p$ and $j >p$ are responsible for the two zero blocks in the above expression of $\mathbf{S}^{(k)}$ while $\mathbf{S}_1$ and $\mathbf{S}_2$ will not be diagonal in general, but only symmetric.} In the second case, we have in general that $\mathcal{C}_{ij} \neq 0$ for $1 \leq i < j$ and therefore for large $k$ we have that $\mathbf{S}^{(k)} \approx \begin{bmatrix}
		\mathbf{D}_1 & \mathbf{0}_{p \times (n-p)} \\
		\mathbf{0}_{(n-p) \times p} & \mathbf{S}_3
	\end{bmatrix}$ where now $\mathbf{D}_1 \in \mathbb{R}^{p \times p}$ is the diagonal with $p$ eigenvalues and we therefore recover their actual $p$ eigenvectors.}
	{\color{black}As previously discussed, we would minimize the objective function of~\eqref{eq:optimization} to zero if and only if $\mathbf{\alpha}$ would contain $p$ eigenvalues of $\mathbf{S}$. But, in general, we cannot assume that these are available a priori.} 
	
	\subsection{The proposed iterative algorithm}
	
In general, the structure in \eqref{eq:theG} is too simple for most applications and therefore we assume now that $\mathbf{\overline{U}}$ in~\eqref{eq:optimization} is written as a product of $k\geq1$ transformations from~\eqref{eq:theG} like:
\begin{equation}
\mathbf{\overline{U}}^{(k)} = \prod_{q=1}^k \mathbf{G}_{i_q j_q} = \mathbf{G}_{i_1 j_1}\mathbf{G}_{i_2 j_2} \dots \mathbf{G}_{i_k j_k},
\label{eq:theU}
\end{equation}
then we can solve~\eqref{eq:theoptimizationproblem} several times: at step $q$ we minimize $\| \mathbf{A} - \mathbf{G}_{i_q j_q}^T \mathbf{S}^{(q-1)} \mathbf{G}_{i_q j_q} \|_F^2$ where $ \mathbf{S}^{(q-1)} = \mathbf{G}_{i_{q-1} j_{q-1}}^T \dots \mathbf{G}_{i_1 j_1}^T \mathbf{S} \mathbf{G}_{i_1 j_1} \dots \mathbf{G}_{i_{q-1} j_{q-1}}$. Similarly to the Jacobi method, our updates are $\mathbf{S}^{(q)} \leftarrow \mathbf{G}_{i_q j_q}^T \mathbf{S}^{(q-1)} \mathbf{G}_{i_q j_q} $.

Based on these results, we give the full description of the proposed method in Algorithm 1. {\color{black}We first compute all the scores $\mathcal{C}_{ij}$. There are $np - \frac{1}{2}p(p-1)$ scores in total (we can think about these scores as being stored in a matrix $\mathcal{C}$ of size $p \times n$ but having the lower triangular part including the main diagonal filled with zero values). Therefore, to compute all these scores takes $9np$ operations -- for all $np$ scores we need to evaluate~\eqref{eq:badcij} and this takes 9 arithmetic operations.

\begin{algorithm}[t]
	\caption{\hspace{-4pt}.\textbf{ Input: }Symmetric matrix $\mathbf{S} \in \mathbb{R}^{n \times n}$, dimension of eigenspace $p \in \mathbb{N}$ such that $p \leq n$, target $\alpha \in \mathbb{R}^p$, and number of basic transformations/iterations $k \in \mathbb{N}$.\newline \textbf{Output: } First $p$ columns of the approximate eigenspace $\mathbf{\overline{U}} \in \mathbb{R}^{n \times n}$ as~\eqref{eq:theU}.}
	\begin{algorithmic}
		\State \textbf{1. }{\color{black}Initialize: $\mathbf{\overline{U}}^{(0)} \leftarrow \mathbf{I}_n$, $\mathbf{S}^{(0)} \leftarrow \mathbf{S}$, compute $\mathcal{C}_{ij}$ for all $i$ and $j$ using~\eqref{eq:badcij} for $1 \leq i \leq p$ and $1\leq i < j \leq n$.}
		
		\State \textbf{2. }Iterative process, for $q = 1,\dots,k:$
		\begin{itemize}
			\item {\color{black}Find the pivot indices $(i_q, j_q) = \underset{i, j;\ i < j \leq n,\ 1\leq i \leq p}{\arg \max}\ \mathcal{C}_{ij}$.}
			\item {\color{black}Build the locally optimum similarity transform $\mathbf{G}_{i_q j_q}$ according to~\eqref{eq:theG} and~\eqref{eq:theGij} using $\mathbf{S}^{(q-1)}$}.
			\item {\color{black}Update the matrix $\mathbf{S}^{(q)} \leftarrow \mathbf{G}_{i_q j_q}^T \mathbf{S}^{(q-1)} \mathbf{G}_{i_q j_q}$ and the approximate eigenspace $\mathbf{\overline{U}}^{(q)} \leftarrow \mathbf{\overline{U}}^{(q-1)}\mathbf{G}_{i_q j_q}$.
				\item Update the row scores $\mathcal{C}_{i_q s}$ for $i_q < s \leq n$ and $\mathcal{C}_{j_q s}$ for $j_q < s \leq n$ if $j_q \leq p$. Then update the column scores $\mathcal{C}_{r i_q}$ for $1 \leq r < i_q $ and $\mathcal{C}_{r j_q}$ for $1 \leq r \leq \min\{j_q, p\}$.}
		\end{itemize}
	\end{algorithmic}
\end{algorithm}
{\color{black}Then, we proceed with the iterative phase for $k$ steps. First, we need to find the pivot indices $(i_q, j_q)$ by searching all the scores $\mathcal{C}$ and this operation, implemented naively, takes $np$ operations. Once we have the pivot indices, we compute the eigenvalue decomposition of size $2 \times 2$ to compute $\mathbf{G}_{i_q j_q}$ and we consider that this takes a constant number of operations $C$. Then we use $\mathbf{G}_{i_q j_q}$ in the similarity transformation to update the approximate eigenspace and these operations take $6n + 6n = 12n$ plus a small constant number of operations. 
	Finally, now that we have the new matrix $\mathbf{S}^{(q)}$ we need to update the scores in $\mathcal{C}$ which takes around $18(n-1)$ operations.}

{\color{black}Therefore, Algorithm 1 performs $9np + k(np + C + 12n + 18n) = 9np + k(np + 30n +C)$ operations. The largest term in this expression is $knp$ and it is caused by the search of the pivot indices in the matrix $\mathcal{C}$. This is a well known computational issues as it is also present in the classic Jacobi method. 
	Storing row maxima (at the extra cost of a new vector of size $2p$), searching for the maximum value can be reduced overall to $k(n+p)$.\footnote{In our implementation, we use the $\max$ function to find the maximum entry in $\mathcal{C}$ as this provided the fastest results in terms of running time.} Therefore, the computational complexity of Algorithm 1, with a careful implementation, is $9np + k(n+p + 30n +C) = 9np + k(31n + p + C)$ operations. As such, Algorithm 1 exhibits complexity $O(nk + np + kp)$. In general, the dominant factor is $O(nk)$, say for example when $k$ is of order $O(n \log n)$, leading to an overall dominant complexity factor of $O(n^2 \log n)$.}
	
	{\color{black}The parameters $n,p$ and $k$ depend on the practical application at hand. In this paper, we usually have that $p \ll n$. The number of iterations/transformations $k$ is very important as it balances, in general, the sparsity of the approximated eigenspace with its accuracy. In this paper, we suggest taking at most $k \sim O(n \log_2 n)$ to balance the accuracy of the eigenspace with its sparsity, see the Experimental Results section. We would like to note that, for $k \sim O(np)$ we expect to have dense approximations of the true eigenspaces and therefore sparsity is lost. The sparsity level and pattern of the resulting transform depend heavily on the indices chosen in the iterative process.

While the proposed method can be used to compute eigenspaces to high accuracy, we note that in this paper, we do not aim to compare against well-established algorithms from the literature that perform eigenvalue decompositions of symmetric matrices. These methods fully exploit block matrix operations, the memory hierarchy and parallelism and therefore, while our proposed method may compare favorably in terms of the number of operations, we expect them to have better (lower) running times. Our method can benefit from a low-level, multi-core careful implementation to significantly reduce running time. Some remarks are in order.}


Similarly to the {\color{black}classical} Jacobi method, Algorithm 1 also has two steps: find pivot indices $(i,j)$ and apply a similarity transformation such that we are closer to the stated objective ({\color{black}diagonalize or block diagonalize} the matrix). We differ from the Jacobi method in both steps: the selection of the indices is made with~\eqref{eq:badcij} instead of the maximum absolute value off-diagonal element and {\color{black}the transformation is not only a Jacobi rotation but it may also be a reflection}. {\color{black}Regarding these scores, for $\mathbf{\alpha} = \mathbf{1}_{p \times 1}$ we have that $\mathcal{C}_{ij} = |S_{ij}|$ like in the Jacobi method when $S_{ii} = S_{jj}$ and $1 \leq i \leq p < j \leq n$.} Despite these differences, our method can also benefit from the parallelization techniques developed for the Jacobi method~\cite{Shroff1990APA,9056494}.

\noindent \textbf{Remark 1 (A block version of Algorithm 1).} In Algorithm 1 we have chosen indices two at a time. We can extend the method to deal with blocks of size $b \times b$ just as in block Jacobi methods~\cite{BECKA2002243,doi:10.1137/090748548,YusakuYamamoto2014}. In this case, we do assume we have a procedure to compute the EVD of a $b \times b$ matrix (the simple formulas of the $2 \times 2$ case are no longer available). {\color{black}For the proposed method, a natural block size is $b = 2p$ such that we pick $p$ pairs of two indices like $(i,\ \underset{j}{\arg \max}\ \mathcal{C}_{ij} )$ for all $i \leq p$ restricting such that no duplicate indices appear. In effect, this means that we pair all allowed indices $i=1,\dots,p$ with a different index $j > p$. This way of choosing the indices is equivalent to setting $\mathbf{\alpha} = \mathbf{1}_{p \times 1}$.}$\hfill \blacksquare$

Given the proposed algorithm, we now provide some theoretical insights into its behavior and show numerical experimental results to validate the approach.

\subsection{Convergence analysis of the proposed algorithm}

Because at each step of the proposed algorithm, we make choices to reduce the objective function in~\eqref{eq:theoptimizationproblem} maximally we are guaranteed to converge. As long as some $\mathcal{C}_{ij} > 0$ progress is possible and these scores~\eqref{eq:badcij} are non-zero at least as long as there is one off-diagonal element $|S_{ij}| > 0$. As the Jacobi method cancels at each step the largest off-diagonal element, its analysis was done differently, using the off-diagonal ``norm'' $\sqrt{\sum_{i=1}^n \sum_{j > i}^n |S_{ij}|^2}$ ($i$ would go only until $p$ and not $n$).

When the given $\mathbf{S}$ is a large sparse matrix we are concerned with the fill-in that happens during Algorithm 1. We note that, at each step, the fill-in is at most $O(n)$. Regarding the scores, we have that most $S_{ij} = 0$ and therefore the critical quantity in~\eqref{eq:badcij} which is $|S_{ii} - S_{jj}| - (S_{ii} - S_{jj})$ takes either the value 0 when the diagonal entries are ordered decreasing, i.e., $S_{ii} > S_{jj}$), or $2(S_{jj} - S_{ii})$ otherwise (in this case, the $\mathbf{G}_{ij}$ in~\eqref{eq:theGij} is a reflector that flips rows/columns $i$ and $j$ in $\mathbf{S}$, ensuring a diagonal with decreasing elements). Therefore, by ordering the diagonal entries of $\mathbf{S}$ we can have exactly as many non-zero scores $\mathcal{C}_{ij}$ as non-zero off-diagonal entries $S_{ij}$.

Maximizing scores leads to the best possible update at each step of the algorithm and can be exploited to build {\color{black}sparse eigenvectors} (we stop after just a few transformations $k \ll n^2$). Still, this step can be expensive as the maximum score needs to be found at each iteration. Furthermore, if we are looking to approximate the eigenspace giving up sparsity, then sweeping all the indices in order is an appropriate solution (analogously to cyclic Jacobi methods~\cite{10.2307/1993275,10.2307/2098943}).


As with the Jacobi method for complex-valued matrices~\cite{10.2307/2098731}, the proposed algorithm is immediately extensible to the Hermitian case. In~\eqref{eq:badcij}, for Hermitian $\mathbf{S} \in \mathbb{C}^{n \times n}$ we have scores $\mathcal{C}_{ij}$ similar to \eqref{eq:badcij} except for the quantity $\sqrt{ (S_{ii} - S_{jj})^2 + 4S_{ij}^2 }$ which becomes $\sqrt{ (S_{ii} - S_{jj})^2 + 4|S_{ij}|^2 }$ and then the decomposition in~\eqref{eq:theGij} is over the complex field.

For our algorithm, which is based on error quantities different from the classic Jacobi algorithm, we can provide the following two results.

{\color{black}\noindent \textbf{Proposition 3.} Assuming $\mathbf{S}$ is a symmetric positive definite matrix whose diagonal entries are ordered decreasing and with eigenvalues $\lambda_1 \geq \lambda_2 \geq \dots \geq \lambda_n > 0$, let $\mathbf{\alpha} = \mathbf{1}_{p \times 1}, p \geq 1$, with condition number $\kappa = \frac{\lambda_1}{\lambda_n} $, then for all pairs of indices $(i,j)$ for which the off-diagonal element obeys $S_{ij}^2 \geq \lambda_n^2$ we have the lower bound:
\begin{equation}
	\mathcal{C}_{ij} \geq \delta \frac{\lambda_n}{\kappa},
\end{equation}
where $\delta \approx 0.6$ is a constant.$\hfill \blacksquare$

Bounding from below the off-diagonal elements $S_{ij}$ is essential to guarantee $\mathcal{C}_{ij} > 0$ and make non-trivial progress in the iterations of the proposed algorithm. Thus, Proposition 3 ensures the minimum progress achieved at each step of the proposed algorithm.

For the overall algorithm, we can provide the following theoretical result.

{\color{black}\noindent \textbf{Proposition 4.} Assuming $\mathbf{S}$ is a symmetric positive definite matrix whose diagonal entries are ordered decreasing and with eigenvalues $\lambda_1 \geq \lambda_2 \geq \dots \geq \lambda_n > 0,$ let $\mathbf{\alpha} = \mathbf{1}_{p \times 1}, p \geq 1$, and assuming that at each step $q$ we choose a pair of indices $(i,j)$ such that the off-diagonal element obeys $(S_{ij}^{(q)})^2 \geq \lambda_n^2$ throughout applying Algorithm 1 which constructs in $k$ steps the eigenspace approximation $\overline{\mathbf{U}}^{(k)}$, then we have:
	\begin{equation}
		\left(  1 +  \frac{\delta}{p n \kappa^2 } \right)^k \text{tr}(\mathbf{AS}) \leq \text{tr}( \mathbf{A}(\overline{\mathbf{U}}^{(k)})^T \mathbf{S} \overline{\mathbf{U}}^{(k)}).
	\end{equation}
	Note that the right-hand side is bounded by $\sum_{t = 1}^p \lambda_t$. The result holds for $k \leq  \left\lfloor \log\left(\frac{\sum_{t=1}^{p} \lambda_t}{\text{tr}(\mathbf{AS})}\right) \log^{-1}\left(1 + \frac{\delta}{pn\kappa^2}\right) \right\rfloor$ or for the smallest $k$ such that $(S_{ij}^{(k)})^2 \leq \lambda_n^2$ for all off-diagonal pairs $(i,j)$. As before, $\delta \approx 0.6$ is a constant.$\hfill \blacksquare$
	
	The condition on the off-diagonal element is essential to establish these theoretical results. The condition could also read $S_{ij}^2 \geq \gamma \lambda_n^2$ for a fixed arbitrary small $\gamma > 0$ and the analysis follows exactly the same steps with this multiplicative factor.
	
	Under the assumptions of Propositions 3 and 4, we show that the proposed method has a linear convergence rate. The factors $n$ and $p$ are natural in the convergence rate because increasing the dimension of the problem or the number of principal components to approximate will lead to an increase in the number of basic transformations needed to reach the same accuracy. To the best of our knowledge, this is the first result that establishes a convergence result in the trace quantity which also depends on the condition number of the input positive symmetric matrix $\mathbf{S}$.
	
	

	
	\subsection{The sparsity of the constructed eigenspace}
	
	The sparsity of the eigenspace (the number of non-zeros and their pattern) heavily depends on the choice of indices $(i,j)$ made as the proposed greedy algorithm progresses. As such, it is impossible to understand in general the number/distribution of the non-zeros. We can provide the following remarks.
	
	\noindent \textbf{Proposition 5.} Starting from a diagonal matrix $\mathbf{D}^{(0)}$ of size $n \times n$, following updates of the form $\mathbf{D}^{(q)} \leftarrow \mathbf{G}_{i_q j_q}^T \mathbf{D}^{(q-1)} \mathbf{G}_{i_q j_q}$ with uniformly random indices $(i_q, j_q)$, assuming no exact cancellations take place in the products, then after $k$ such steps the expected number of non-zeros on the off-diagonal approximately follows a sigmoid-shaped curve:
	\begin{equation}
		\frac{ O(n^2) }{ 1 + \exp(- k / O(n) ) }.\hfill \blacksquare
	\end{equation}

	\noindent \textbf{Remark 2 (Worst-case fill-in).} In the worst-case fill-in scenario, just $k = n-1$ suffice to construct an eigenspace with no zero entries. This is achieved for the choices $i_q = 1,\dots, n-1$ and $j_q = i_q +1$ at step $q$.$\hfill \blacksquare$
	
	These two results suggest that, on average, we expect slow fill-in for low values of $k$ in Algorithm 1. Still, note that due to the proposed factorization we can construct numerically simple eigenspaces, i.e., low $k$, but that do not exhibit low sparsity.

	\section{Experimental results}
	
	In this section, we use the proposed method to construct low-rank approximations of 1) random matrices; 2) Laplacian matrices of random graphs; 3) covariance matrices in the context of Principal Component Analysis (PCA) as a dimensionality reduction technique before classification. {\color{black}For a given symmetric $\mathbf{S}$, to measure the quality of the eigenspace approximation $\mathbf{\overline{U}}$ we use the measure:}
	\begin{equation}
		{\color{black}\epsilon(\mathbf{S}, \mathbf{\overline{U}}) = \frac{ \text{tr}( |\mathbf{\overline{U}}^T \mathbf{S\overline{U}} | ) }{\sum_{t=1}^p |\lambda_t|}.}
		\label{eq:error}
	\end{equation}
	The eigenvalues $\lambda_t$ can be either the lowest or highest of $\mathbf{S}$ and we assume $\sum_{t=1}^p |\lambda_{t}| \neq 0$, i.e., we are not trying to recover the null space of $\mathbf{S}$ of size $p$. Good approximations are obtained whenever $\epsilon \approx 1$ (or $100\%$). {\color{black}We will always compare our proposed method against the classic Jacobi method (the pivot element has the maximum absolute value among the off-diagonal elements and besides rotations we also use reflections, see~\eqref{eq:theGij}).} Full source code for the proposed method is available online\footnote{https://github.com/cristian-rusu-research/JACOBI-PCA}.
	\begin{figure}[!t]
		\centering
		\includegraphics[trim = 45 20 51 30, clip, width=0.5\textwidth]{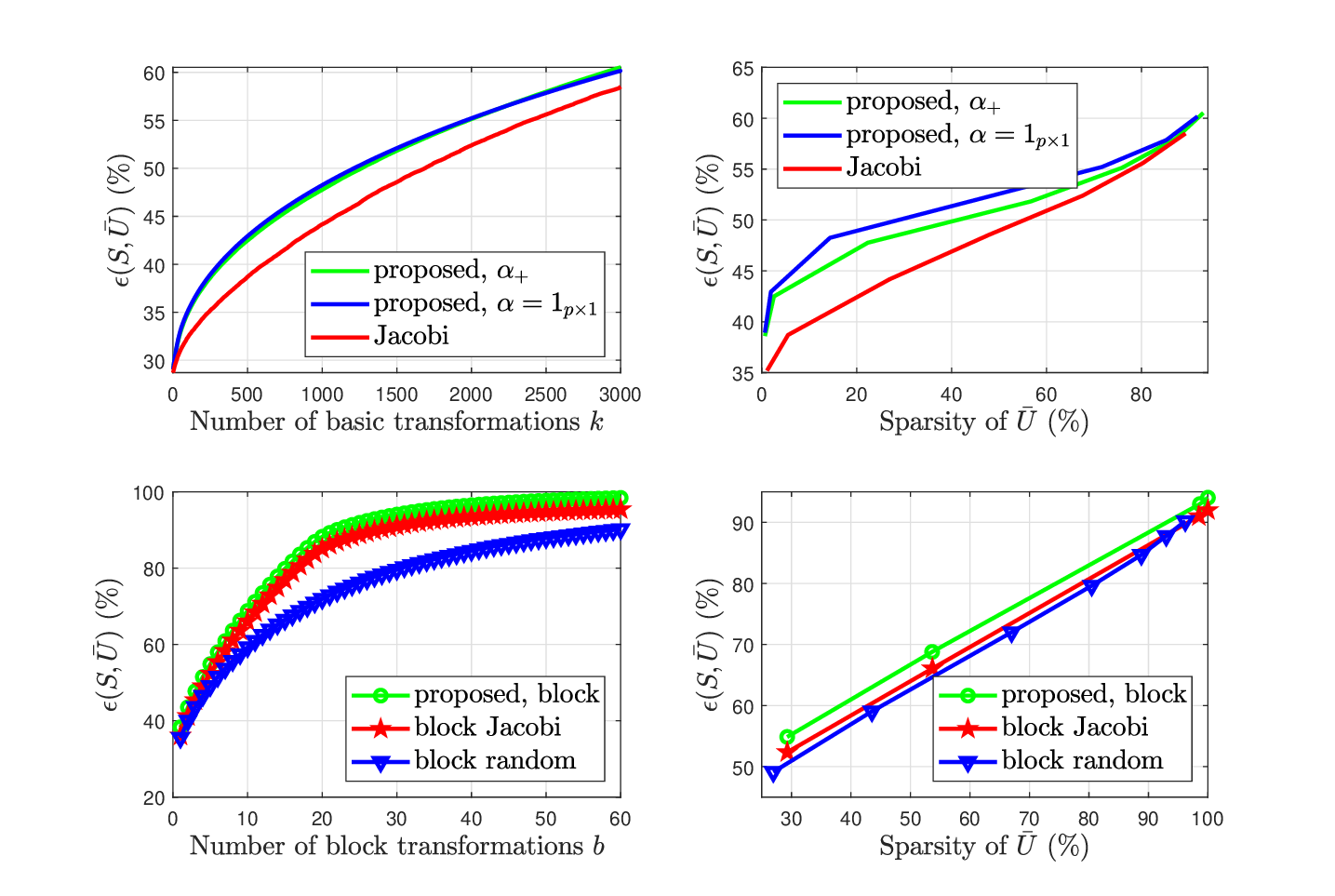}
		\caption{Average approximation accuracy~\eqref{eq:error} over 10 realizations for random symmetric matrices of size $n = 1024$ when recovering $p = 20$ eigenvectors associated with the highest eigenvalues. For the proposed method we show simulations for $\mathbf{\alpha} = \mathbf{1}_{p \times 1}$ and $\mathbf{\alpha}_{+} = \log_2(p+1,\dots,2)$ (a positive decreasing sequence). Figures show in order: top left the evolution of~\eqref{eq:error} with the number of transformations $k$ in the approximation $\mathbf{\overline{U}}$; top right the sparsity of the eigenspace ({\color{black}number of non-zero elements in $\mathbf{\overline{U}}$ divided by the total number of elements which is $np$}); bottom left the evolution of~\eqref{eq:error} for the block variant of the proposed method with block size $b  = 2p$; bottom right the sparsity of the eigenspace created via the block proposed method (see Remark 1). {\color{black}The block random method refers to performing the best $2p\times 2p$ similarity transformations on randomly chosen pivot indices.}}
		\label{fig:figure1}
	\end{figure}
	
	\subsection{Low-rank approximations of random symmetric matrices}
	
	{\color{black}We generate random matrices $\mathbf{X} \in \mathbb{R}^{n \times n}$ with entries drawn independent and identically distributed from the standard Gaussian distribution, then initialize $\mathbf{S} \leftarrow \mathbf{X} \mathbf{X}^T$ and we apply Algorithm 1 with $p=20$ and various $k$ on $\mathbf{S}$.} {\color{black}In Figure~\ref{fig:figure1} we show average results obtained with Algorithm 1 on random symmetric matrices for $n = 1024$. We compare against the Jacobi method and also show the block variant of the proposed method (see Remark 1).} {\color{black}In Figure~\ref{fig:figure1_new} we perform the same numerical experiments as in Figure~\ref{fig:figure1} with the distinction that now the matrix $\mathbf{S}$ has rank $p$ (use the truncated singular value decomposition to make $\mathbf{S}$ low rank before applying Algorithm 1).} We always perform better than the Jacobi method in these experiments, showing that selecting pivot indices via~\eqref{eq:badcij} brings benefits.
	
	An important point of Figure~\ref{fig:figure1} is to show the accuracy of the approximating eigenspaces as a function of the sparsity of these eigenspaces. High sparsity leads, in general, to poorer results highlighting the fact that the true eigenspaces are indeed not sparse and therefore no high-quality approximation can be computed for low $k$.
	
	{\color{black}
		We have chosen the trace term~\eqref{eq:error} because this is the quantity that we have seen in the theoretical results described in this paper. There are, of course, several other ways to estimate the quality of the experimental results. Given the true eigenspace $\mathbf{U}_p$ of size $n \times p$ and the accuracy achieved by Algorithm 1 $\mathbf{\overline{U}}_p$ also of size $n \times p$ we can define the maximum canonical angle between subspaces spanned by the columns of $\mathbf{U}_p$ and $\mathbf{\overline{U}}_p$~\cite{10.2307/2005662}. We simply call this the angle between subspaces $\mathbf{U}_p$ and $\mathbf{\overline{U}}_p$:
		\begin{eqnarray}
			\text{angle}(\mathbf{U}_p, \mathbf{\overline{U}}_p) = \sin^{-1} ( \| \mathbf{U}_p - \mathbf{\overline{U}}_p(\mathbf{\overline{U}}_p^T \mathbf{U}_p)  \|_2 ).
			\label{eq:subspacedistance}
		\end{eqnarray}
		We have that $0 \leq \text{angle}(\mathbf{U}_p, \mathbf{\overline{U}}_p) \leq \pi/2$. When the distance is $\pi/2$ then the two subspaces are orthogonal to each other and at zero the two subspaces are the same.
		
		In Figures~\ref{fig:figure_full_rank} and~\ref{fig:figure_low_rank} we show the accuracy~\eqref{eq:error} and subspace distance~\eqref{eq:subspacedistance} for random symmetric positive semidefinite matrices of full rank and of rank $p$, respectively. Unlike the previous experimental setup, this time we allow large $k$ and $b$, i.e., we forsake sparsity for the pursuit of accuracy. This situation is akin to the sweeps, i.e., $O(n^2)$ similarity transformations, done in the classic Jacobi method. {\color{black}We note that Algorithm 1 performs much better when the input matrices are themselves low rank (see Figure~\ref{fig:figure_full_rank} versus Figure~\ref{fig:figure_low_rank}). Differences between the methods are especially significant in the subspace distance.} Furthermore, we observe that the two performance indicators are correlated: higher trace quantities lead to lower subspace distances. But, we do want to point out that for any significant reduction in the subspace distance we need high values of the approximation accuracy: only $\epsilon \approx 1$ leads to consistently small angles in~\eqref{eq:subspacedistance} (see the right side versus the left side in Figures~\ref{fig:figure_full_rank} and~\ref{fig:figure_low_rank}). {\color{black}As we expect that in general the true eigenspaces are not actually sparse, we need a large number of transformations $k$ to achieve good approximation results. In these cases, sparsity is lost almost entirely.} For large $b$, the block variant of Algorithm 1 performs well in these scenarios.}
	\begin{figure}[!t]
		\centering
		\includegraphics[trim = 45 15 51 30, clip, width=0.5\textwidth]{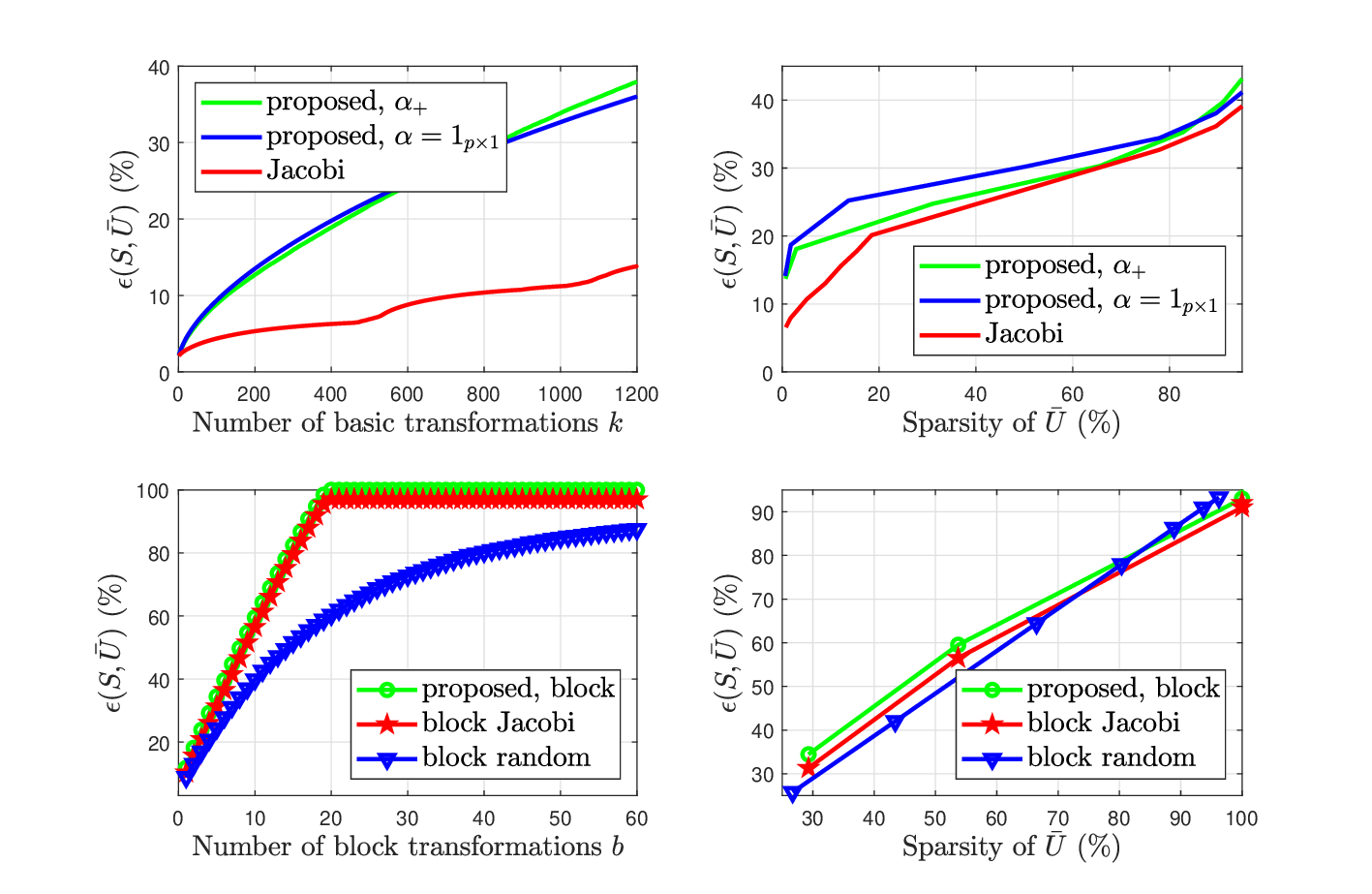}
		\caption{Results for the same experimental setup as in Figure~\ref{fig:figure1} but now the random symmetric matrices which are the input to Algorithm 1 have rank $p$.}
		\label{fig:figure1_new}
	\end{figure}
	\begin{figure}[!t]
		\centering
		\includegraphics[trim = 45 15 55 15, clip, width=0.5\textwidth]{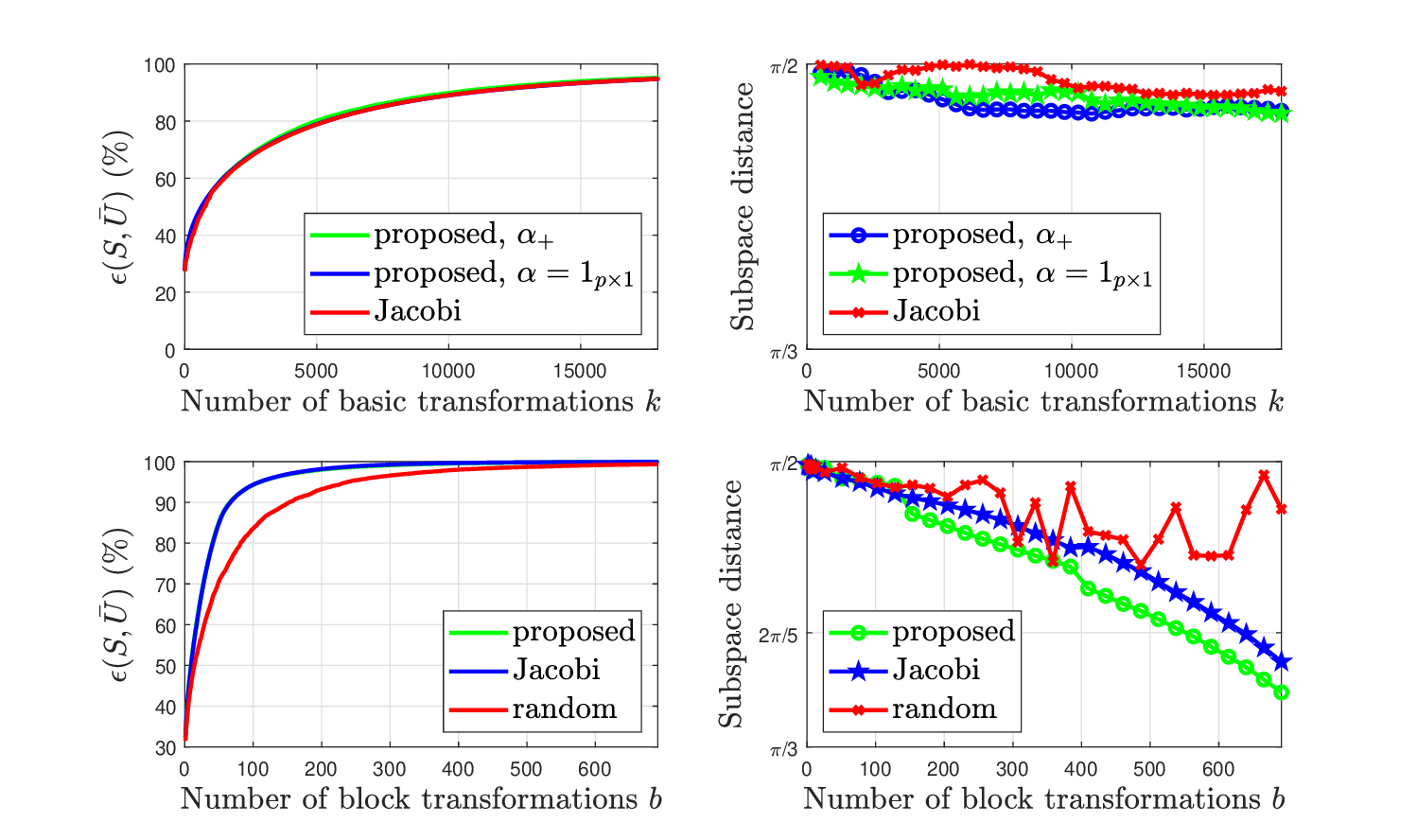}
		\caption{For the same experimental setup as in Figure~\ref{fig:figure1} (random symmetric positive semidefinite matrices) we show the evolution of the average approximation accuracy~\eqref{eq:error} and the subspace distance~\eqref{eq:subspacedistance} for a large number of basic transformations $k$ (top row) and a large number of block transformations $b$ (bottom row). The random method refers to choosing the pivot indices $(i,j)$ randomly.}
		\label{fig:figure_full_rank}
	\end{figure}
	\begin{figure}[!t]
		\centering
		\includegraphics[trim = 70 20 65 30, clip, width=0.5\textwidth]{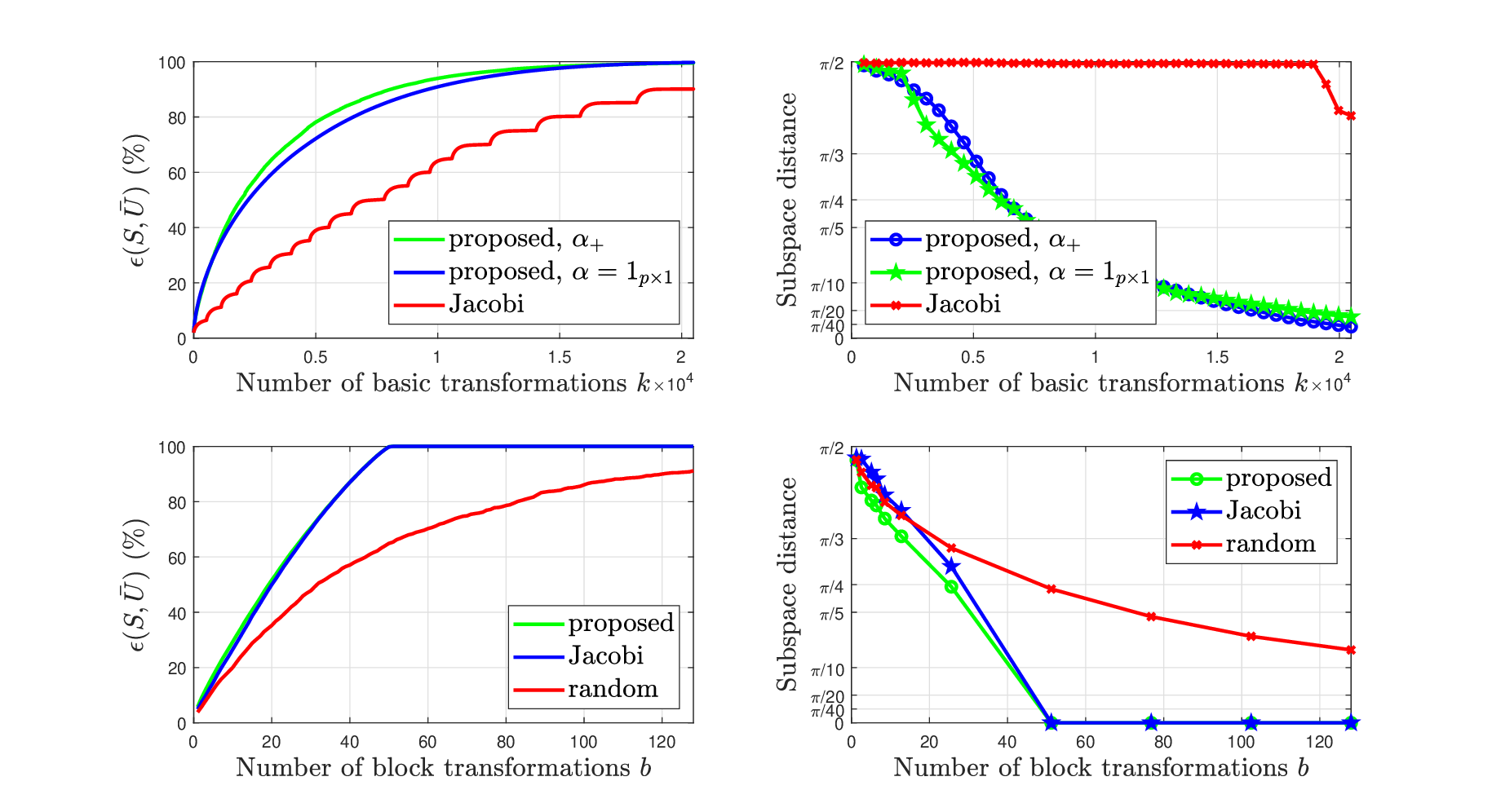}
		\caption{For the same experimental setup as in Figure~\ref{fig:figure1_new} (low-rank random symmetric positive semidefinite matrices) we show the evolution of the average approximation accuracy~\eqref{eq:error} and the subspace distance~\eqref{eq:subspacedistance} for a large number of basic transformations $k$ (top row) and a large number of block transformations $b$ (bottom row).}
		\label{fig:figure_low_rank}
	\end{figure}
	
	{\color{black}Finally, in Figure~\ref{fig:figure5} we also show the individual convergence to the largest $p$ eigenvalues for random symmetric matrices. Note how in this case, larger eigenvalues are recovered better, but results will depend in general on the choice of $\mathbf{\alpha}$.}
	

	\subsection{Graph Fourier transforms: eigenvectors of Laplacians}
	
	{\color{black}In many graph signal processing applications~\cite{Sandryhaila20131644,8347162,STANKOVIC2020102802,DALCOL2022108690} eigenvalue decompositions play an important role.
		
		We define an undirected graph $\mathcal{G}$ as the pair $(\mathcal{V}, \mathcal{E})$ where $\mathcal{V}$ is the set of size $n$ representing the graph nodes and $\mathcal{E}$ is the set of size $m$ denoting the graph edges. Based on these two sets we can define: the $n \times n$ diagonal matrix $\mathbf{D}$ called the degree matrix (the $i\text{th}$ element on the diagonal is the degree of node $i$, i.e., the number of edges connected to it) and the $n \times n$ adjacency matrix $\mathbf{A}$ whose elements are zero except for entries $(i,j)$ and $(j,i)$ which have the value one if and only if an edge exists between nodes $i$ and $j$. In general, we assume that $k \ll n^2$ and therefore $\mathbf{A}$ is a sparse matrix. Based on these two matrices we can define the symmetric positive semidefinite Laplacian matrix $\mathbf{L} = \mathbf{D} - \mathbf{A}$. In this section, we are interested in estimating the eigenvectors associated with the lowest eigenvalues of an undirected graph Laplacian $\mathbf{L}$ as they contain important information about the structure of the graph $\mathcal{G}$.

		The eigenvectors of the Laplacian $\mathbf{L}$ define the graph Fourier transform~\cite{RICAUD2019474} and are useful as they provide ``low-frequency'' information about graph signals, in the context of graph signal processing.} We generate random community graphs of $n = 256$ nodes with the Graph Signal Processing Toolbox\footnote{https://epfl-lts2.github.io/gspbox-html/} for which we decompose the sparse positive semidefinite Laplacians as $\mathbf{L} = \mathbf{U} \text{diag}(\mathbf{\lambda}) \mathbf{U}^T$ and recover the eigenvectors from $\mathbf{U}$ associated with the lowest eigenvalues from $\text{diag}(\mathbf{\lambda})$. {\color{black}To compute the eigenspace of the lowest eigenvalues we take Algorithm 1 with $\mathbf{\alpha} = -\mathbf{1}_{p \times 1}$ and $\mathbf{\alpha}_{-}$ a negative increasing sequence. This is equivalent to applying Algorithm 1 on the input matrix $-\mathbf{L}$, instead of $\mathbf{L}$. {\color{black}Results on the Laplacian matrix $\mathbf{L}$ are shown in Figure~\ref{fig:figure2}. As we approximate the eigenvectors associated with the lowest positive eigenvalues, the first $p$ diagonal entries of $\mathbf{S}^{(k)}$ will decrease as $k$ increases, and we therefore expect to monotonically decrease towards $\epsilon = 100\%$ in \eqref{eq:error}.}}
	\begin{figure}[!t]
		\centering
		\includegraphics[trim = 75 0 65 15, clip, width=0.5\textwidth]{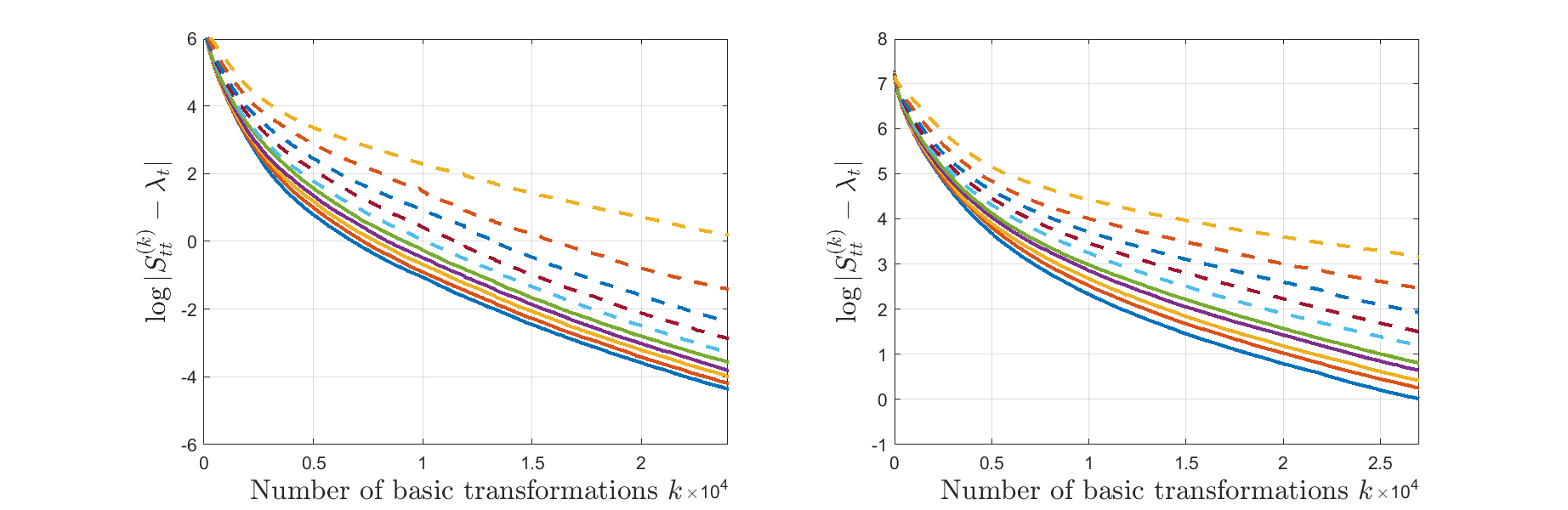}
		\caption{Average distance between the diagonal elements $S_{tt}^{(k)}$ and the true eigenvalues $\lambda_t$, with $t=1,\dots,p$, averaged over 10 realizations of full rank random symmetric positive semidefinite matrices of size $n = 128$ (left) and $n = 256$ (right). We have $p = 10$ and $\mathbf{\alpha}_{+} = \log_2(p+1,\dots,2)$. Solid lines show the distances to the larger 5 eigenvalues and the dashed lines are the distances for the other, lower, 5 eigenvalues.}
		\label{fig:figure5}
	\end{figure}
	\begin{figure}[!t]
		\centering
		\includegraphics[trim = 50 0 55 15, clip, width=0.5\textwidth]{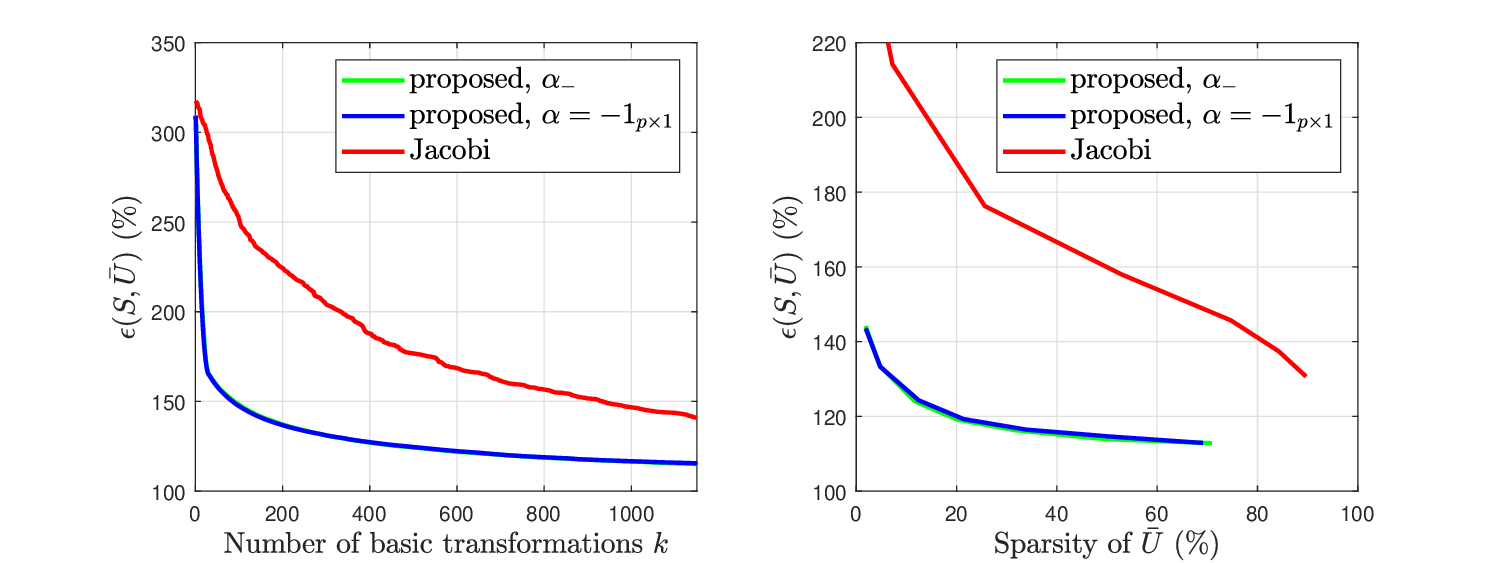}
		\caption{Average approximation accuracy~\eqref{eq:error} over 10 realizations for random community graphs of size $n = 256$ created via the GSP Toolbox when recovering $p = 32$ eigenvectors associated with the lowest eigenvalues of the graph Laplacians $\mathbf{L}$. We show accuracy as a function of the number of transformations $k$ (left) and the sparsity of the eigenspace (right). Here, lower is better and $100\%$ means perfect recovery of all $p$ smallest eigenvalues. We show simulations for $\mathbf{\alpha} = -\mathbf{1}_{p \times 1}$ and $\mathbf{\alpha}_{-} = -\log_2(p+1,\dots,2)$ (a negative increasing sequence).}
		\label{fig:figure2}
	\end{figure}

	\subsection{Dimensionality reduction via sparse PCA}
	
	{\color{black}
		In the context of machine learning applications, given a dataset, it is of interest to compute the eigenvectors associated with the highest eigenvalues of the covariance matrix -- a task called Principal Component Analysis (PCA). Recently, random projections have been widely used to perform dimensionality reduction on large-dimensional datasets successfully. This is therefore a good candidate scenario where Algorithm 1 can be applied to perform dimensionality reduction using sparse projections, i.e., computing $\mathbf{\overline{U}}_p$ by Algorithm 1 with low $k$.}
	
	We consider a classification example with the USPS dataset\footnote{https://github.com/darshanbagul/USPS\_Digit\_Classification} with 10 classes for which $n = 256$ and the number of data points is $N = 9298$. From the data matrix $\mathbf{X} \in \mathbb{R}^{n \times N}$ we explicitly compute the covariance $\mathbf{C} = \mathbf{X} \mathbf{X}^T$ and apply Algorithm 1 on $\mathbf{C}$ to get an approximation of the $p$ principal components. We use the $K$-nearest neighbors ($K$-NN) classification algorithm with $K = 15$, but before we perform dimensionality reduction with $p = 15$ to improve the running time of the classifier. Results are shown in Figure~\ref{fig:figure3}. We compare against the full PCA, the Sparse Johnson-Lindenstrauss (S-JL) transform~\cite{10.1145/2559902} for $p = 30$ components (for $p = 20$ as with PCA the performance is poor because S-JL is not data dependent) with sparsity $s = 3$ per column (overall sparsity $30\%$).
	
	\begin{figure}[!t]
		\centering
		\includegraphics[trim = 45 0 43 12, clip, width=0.5\textwidth]{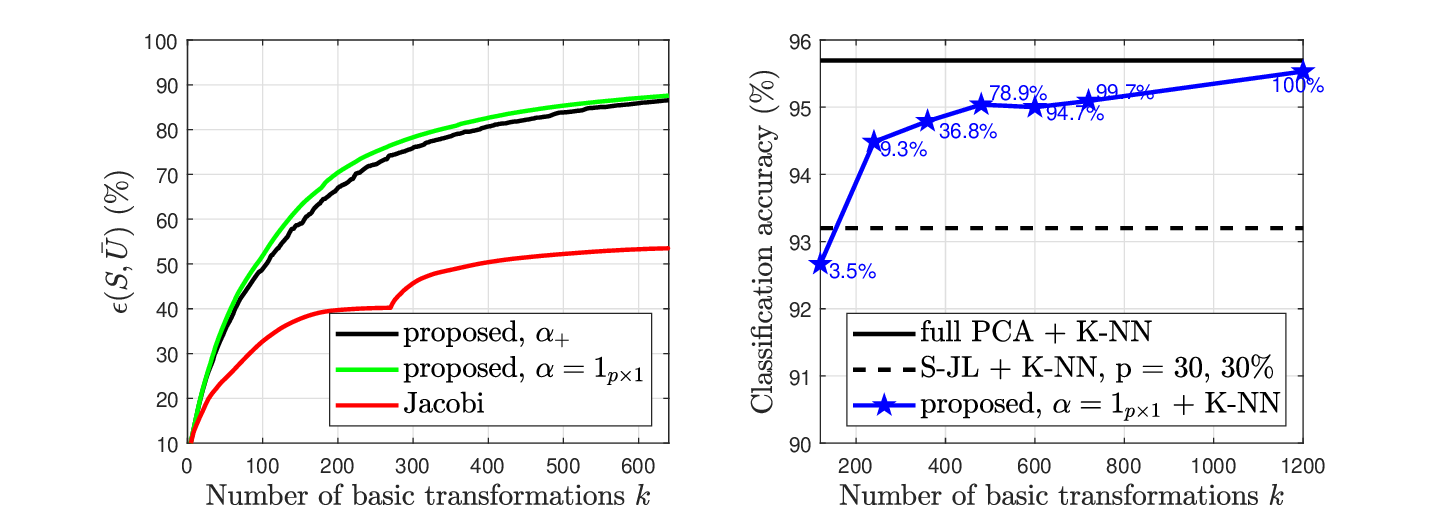}
		\caption{Average approximation accuracy~\eqref{eq:error} and average classification accuracy, on the left and right, respectively, over 10 random realizations obtained by the proposed algorithm and with $K$-NN for the USPS dataset split into $N_\text{train} = 8000$ and $N_\text{test} = 1298$ after dimensionality reduction with $p = 20$ was performed. Percentages on the plot shown in blue represent the sparsity of the eigenspace. {\color{black}For the plot of the right, for the 7 approximations we build the trace quantities $\epsilon$ from \eqref{eq:error} are: 55.1\%, 69.9\%, 75.2\%, 78.4\%, 80.9\%, 82.3\%, and 86.2\%.}}
		\label{fig:figure3}
	\end{figure}
	\begin{figure}[!t]
		\centering
		\includegraphics[trim = 75 0 61 23, clip, width=0.5\textwidth]{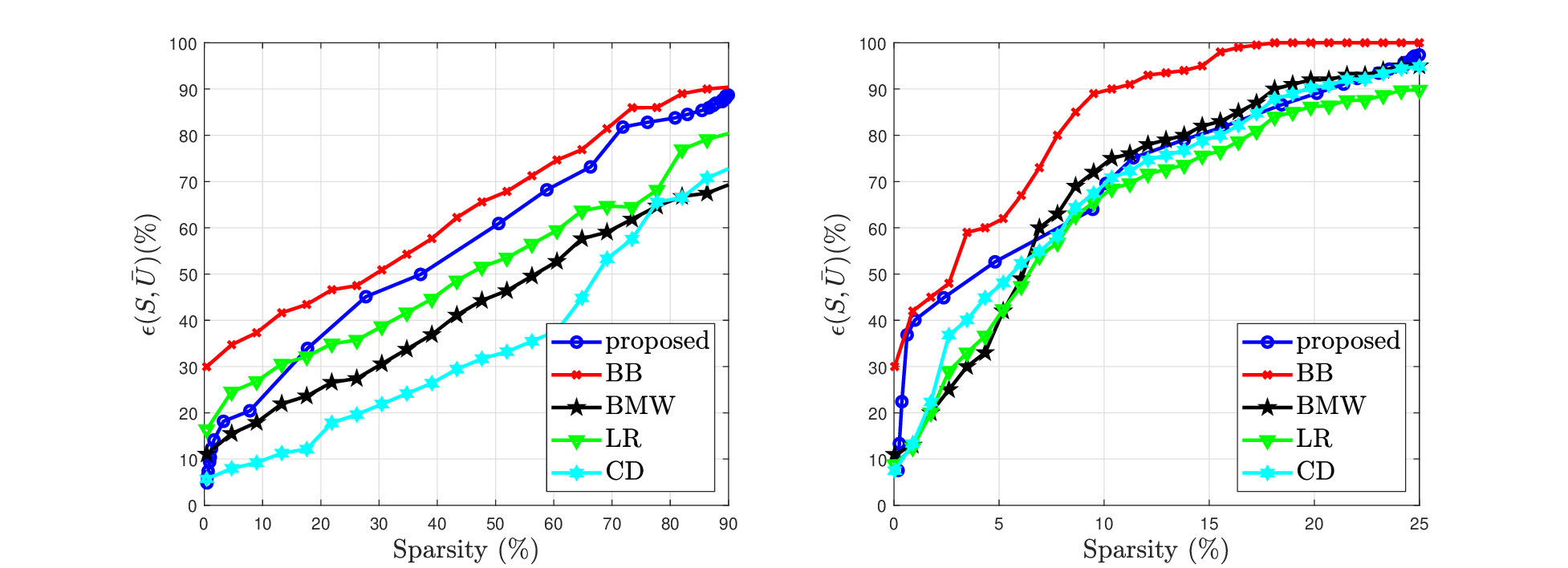}
		\caption{Comparison of sparse PCA methods from the literature for the USPS (left) and MNIST (right) datasets for $p = 20$ principal directions of various sparsity levels. Note how, for all considered methods, accuracy~\eqref{eq:error} depends on the sparsity of the computed eigenspace.}
		\label{fig:sparse_pca}
	\end{figure}
	{\color{black}Finally, in Figure~\ref{fig:sparse_pca} we compare our proposed method against algorithms for sparse PCA recently presented in the literature: branch and bound (BB)~\cite{Berk2019}, bipartite matching (BMW)~\cite{NIPS2015_2b8a6159}, low-rank updates (LR)~\cite{papailiopoulos2013sparse}, and coordinate descent (CD)~\cite{Beck2016}. The BB method performs best but is the slowest, for example, for sparsity level 50\% the accuracy is almost 70\% but the algorithm ran for approximately 15 minutes. Our proposed method performed second best in these tests and in terms of running times, we have 0.3 seconds for the best performing projection with sparsity 100\% and accuracy 85\%, and 0.05 seconds for the projection with sparsity 2\% and accuracy 20\%. We also note that our proposed method does not reach high accuracy for dense eigenspaces (100\% sparsity leads only to 85\% accuracy). This is because numerical experiments show that our proposed method needs a large number of basic transformations to achieve very high accuracy (either~\eqref{eq:error} or~\eqref{eq:subspacedistance}).}

	\section{Conclusions}
	{\color{black}In this paper, we have described an extension of the Jacobi method for the diagonalization of symmetric (or Hermitian) matrices for the computation of just a few extreme approximate eigenvalue/eigenvector pairs.
		
		While the optimization problem we described is non-convex and therefore hard to solve optimally in general, the proposed method is based on new closed-form updates that guarantee locally optimal error reduction. {\color{black}The new approach is possible because we consider $2 \times 2$ orthogonal rotation and reflector transformations together and rely on Procrustes results.} {\color{black}Furthermore, the proposed method uses a relatively small number of computationally simple steps, can be trivially parallelized, and exhibits a low memory footprint.} Furthermore, the theoretical results we provide show the linear convergence of the method. We show experimental results where we recover {\color{black}sparse eigenvectors} associated with both the highest and lowest eigenvalues of symmetric matrices and highlight a general trade-off between sparsity and eigenspace recovery accuracy. {\color{black}The proposed method performs well when considering two performance indicators: trace error terms and principal angles between subspaces.}
		
		Finally, we show applications to graph signal processing and data processing where we compare against the classic Jacobi diagonalization method and compare favorably against sparse principal component analysis methods from the literature.}
	
	
	
	
	
	\appendices
	
	\section*{Appendix: Proof of Proposition 1}
	
	We denote the left matrix by $\mathbf{A}$ and then we have the objective function $\| \mathbf{A} - \mathbf{\overline{U}}^T \mathbf{S} \mathbf{\overline{U}} \|_F^2 = \|\mathbf{A} \|_F^2 + \| \mathbf{S} \|_F^2 - 2\text{tr}(\mathbf{A} \mathbf{\overline{U}}^T \mathbf{S} \mathbf{\overline{U}})$. {\color{black}By \cite[Theorem~19]{Inequalities} and \cite[Theorem~1.2]{doi:10.1137/0719089} (which is based on the Courant-Fischer theorem \cite[Corollary 4.3.39]{Horn}), {\color{black}we maximize the trace term by choosing orthogonal $\mathbf{\overline{U}}$ to have $p$ eigenvectors of $\mathbf{S}$ as its first $p$ columns such that} $\text{tr}(\mathbf{A} \mathbf{\overline{U}}^T \mathbf{S} \mathbf{\overline{U}}) = \sum_{i=1}^p \alpha_i \lambda_{\sigma(i)}$ is largest for some ordering $\sigma$. To maximize the trace, the ordering is such that either the largest eigenvalues are picked (when $\alpha_i$ is positive) or the smallest eigenvalues (when $\alpha_i$ is negative).} {\color{black}This ordering $\sigma$ determines the way the eigenvectors of $\mathbf{S}$ are arranged as the first $p$ columns of $\mathbf{\overline{U}}$}.
	
	\section*{Appendix: Proof of Proposition 2}
	
	We use the structured matrices introduced in~\cite{FastPCA}, the proof outline given for Theorem 1 of~\cite{rusu2020constructing} and develop the objective function in~\eqref{eq:theoptimizationproblem} using the definition of the Frobenius norm to
	{\color{black}
		\begin{equation}
			\begin{aligned}
				\|\mathbf{A} - \mathbf{G}_{i j}^T \mathbf{S} \mathbf{G}_{i j} \|_F^2 = & \|\mathbf{A}\|_F^2 + \| \mathbf{S} \|_F^2 - 2\text{tr}(\mathbf{G}_{i j}^T \mathbf{S} \mathbf{G}_{i j} \mathbf{A}  )\\
				= & \| \mathbf{\alpha} \|_2^2 + \| \mathbf{\lambda} \|_2^2 - 2\text{tr}(\mathbf{Z}) - 2\mathcal{C}_{i j},
			\end{aligned}
	\end{equation}}
	{\color{black}where we have used the fact that the Frobenius norm is unitarily invariant and we have defined $\mathbf{Z} = \mathbf{AS}$ for which $\text{tr}(\mathbf{AS}) = \sum_{t=1}^p \alpha_t S_{tt} $. We focus on the cost}
	\begin{equation}
		\mathcal{C}_{i j} = \text{tr}(\mathbf{\widetilde{G}}_{ij}^T \mathbf{\widetilde{S}}_{ ij } \mathbf{\widetilde{G}}_{ij}\  \text{diag}(\begin{bmatrix}
			A_{ii} &
			A_{jj}
		\end{bmatrix})) - Z_{i i} \! -\! Z_{j j}.
	\end{equation}
	With our choice for $\mathbf{\widetilde{G}}_{ij} = \mathbf{V}$ from~\eqref{eq:theGij} and using the EVD of $ \mathbf{\widetilde{S}}_{ ij }$, which is called the Procrustes solution and was shown to be optimal~\cite{Schnemann1968OnTO}, the trace term develops to
	\begin{equation}
		\begin{aligned}
			& \text{tr}(\mathbf{V}^T \mathbf{\widetilde{S}}_{ ij } \mathbf{V}\ \text{diag}(\begin{bmatrix}
				A_{ii} &
				A_{jj}
			\end{bmatrix}))
			= \text{tr}(\mathbf{V}^T \mathbf{V} \mathbf{D}  \mathbf{V}^T \mathbf{V}\ \text{diag}(\begin{bmatrix}
				A_{ii} &
				A_{jj}
			\end{bmatrix})) \\
			& =\text{tr}(\mathbf{D}\ \text{diag}(\begin{bmatrix}
				A_{ii} &
				A_{jj}
			\end{bmatrix})) \\
			& =  \text{tr}(\text{diag}(\mathbf{d}) \text{diag}(\begin{bmatrix}
				A_{ii} &
				A_{jj}
			\end{bmatrix})) \\
			& = \begin{bmatrix}
				A_{ii} & A_{jj}
			\end{bmatrix} \mathbf{d}.
		\end{aligned}
	\end{equation}
	{\color{black}The diagonal matrix $\mathbf{D} = \text{diag}(\mathbf{d})$ is from~\eqref{eq:theGij} and contains two real eigenvalues of the symmetric matrix $\mathbf{\widetilde{S}}_{ij }$. Therefore, the total {\color{black}score} is
		\begin{equation}
			\begin{aligned}
				\mathcal{C}_{i j} = & \begin{bmatrix}
					A_{ii} & A_{jj}
				\end{bmatrix}  \mathbf{d} - Z_{i i} - Z_{j j}  \\
				= & \begin{bmatrix}
					A_{ii} & A_{jj}
				\end{bmatrix}  \mathbf{d} - \begin{bmatrix}
					A_{ii} & A_{jj}
				\end{bmatrix} \begin{bmatrix}S_{i i} \\ S_{j j} \end{bmatrix}\\
				= & \begin{bmatrix}
					A_{ii} & A_{jj}
				\end{bmatrix} \begin{bmatrix}
					d_1 - S_{ii} \\
					d_2 - S_{jj}
				\end{bmatrix}.
			\end{aligned}
		\end{equation}
		Note that $Z_{i i} = A_{ii} S_{i i} $ and $Z_{j j} = A_{jj} S_{j j}$ and the eigenvalues of $\mathbf{\widetilde{S}}_{ i j}$ are stored in $\mathbf{d}$ and are calculated using explicit formulas: $d_{1,2} \in \{ \frac{1}{2}(S_{ii} + S_{jj} \mp \sqrt{ (S_{ii} - S_{jj})^2 + 4S_{ij}^2  } )\}$. Now we have the question of how to order the eigenvalues in $\mathbf{d}$. {\color{black}To ensure $\mathcal{C}_{ij} \geq 0$ always, we have two cases:
			\begin{itemize}
				\item[a)] if $A_{ii} \leq A_{jj}$ then we have $d_{1,2} = \frac{1}{2}(S_{ii} + S_{jj} \mp \sqrt{ (S_{ii} - S_{jj})^2 + 4S_{ij}^2  } )$ and then:
				\begin{equation}
					\begin{aligned}
						\mathcal{C}_{ij} =  & \begin{bmatrix}
							A_{ii} & A_{jj}
						\end{bmatrix} \begin{bmatrix}
							-\frac{1}{2} ( S_{i i } - S_{ j j } + \sqrt{ (S_{i i } - S_{ j j })^2 + 4S_{ij}^2} ) \\
							\frac{1}{2} ( S_{i i } - S_{ j j } + \sqrt{ (S_{i i } - S_{ j j })^2 + 4S_{ij}^2} )
						\end{bmatrix}  \\
						= & {\color{black} \frac{1}{2}  (A_{jj} - A_{ii}) \left( S_{i i } - S_{ j j } + \sqrt{ (S_{i i } - S_{ j j })^2 + 4S_{ij}^2} \right).}
					\end{aligned}
				\end{equation}
				
				\item[b)] if $A_{ii} > A_{jj}$ then we have $d_{1,2} = \frac{1}{2}(S_{ii} + S_{jj} \pm \sqrt{ (S_{ii} - S_{jj})^2 + 4S_{ij}^2  } )$ and then:
				\begin{equation}
					\begin{aligned}
						\mathcal{C}_{ij} =  & \begin{bmatrix}
							A_{ii} & A_{jj}
						\end{bmatrix} \begin{bmatrix}
							-\frac{1}{2} ( S_{i i } - S_{ j j } - \sqrt{ (S_{i i } - S_{ j j })^2 + 4S_{ij}^2} ) \\
							\frac{1}{2} ( S_{i i } - S_{ j j } - \sqrt{ (S_{i i } - S_{ j j })^2 + 4S_{ij}^2} )
						\end{bmatrix}  \\
						= & {\color{black} \frac{1}{2}  (A_{jj} - A_{ii}) \left( S_{i i } - S_{ j j } - \sqrt{ (S_{i i } - S_{ j j })^2 + 4S_{ij}^2} \right).}
					\end{aligned}
				\end{equation}
			\end{itemize}
		} Therefore, the minimizer of~\eqref{eq:theoptimizationproblem} is given in~\eqref{eq:theGij} with the appropriate ordering of the eigenvalues on the diagonal of $\mathbf{D}$.
		
		\section*{Appendix: Proof of Proposition 3}
		
		We are in the case of $A_{ii} > A_{jj}$ in Proposition 2 and we rewrite the value of the scores as:
		\begin{equation}
			\begin{aligned}
				\mathcal{C}_{ij} = & \frac{1}{2} ( \sqrt{ (S_{ii} - S_{jj} )^2 + 4S_{ij}^2} - (S_{ii} - S_{jj})  ) \\
				= &  \frac{2S_{ij}^2}{   \sqrt{ (S_{ii} - S_{jj} )^2 + 4S_{ij}^2} + (S_{ii} - S_{jj})}.
			\end{aligned}
		\end{equation}
		We will upper bound the denominator by using the following bounds: i) $0 < \lambda_n \leq S_{tt} \leq \lambda_1$ for all $t = 1,\dots,n$ which is true by the Rayleigh quotient bounds applied to the standard basis vectors $\mathbf{e}_t$; ii) $-(\lambda_1 - \lambda_n) \leq S_{ii} - S_{jj} \leq \lambda_1 - \lambda_n$ which is true again by the Rayleigh quotient bounds; and iii) $S_{ij}^2 \leq \lambda_1^2$ which true because of the matrix Cauchy-Schwarz inequality $S_{ij}^2 = ( \mathbf{e}_i^T \mathbf{S}\mathbf{e}_j)^2 \leq (\mathbf{e}_i^T \mathbf{S}\mathbf{e}_i) (\mathbf{e}_j^T \mathbf{S}\mathbf{e}_j)  = S_{ii} S_{jj} \leq \lambda_1^2$. 
		Therefore, the denominator is bounded above by
		\begin{equation}
			\! \sqrt{\! \! (S_{ii} - S_{jj} )^2 \! + 4S_{ij}^2} +\! (S_{ii} \! - \! S_{jj}) \leq \! \! \sqrt{(\lambda_1 \! - \! \lambda_n)^2 + 4\lambda_1^2 } +\! (\lambda_1 \! - \! \lambda_n).
		\end{equation}
		The numerator is bounded below by the assumption $S_{ij}^2 \geq \lambda_n^2$. Combining the two results, we reach the following:
		\begin{equation}
			\begin{aligned}
				\mathcal{C}_{ij} \geq & \frac{2 \lambda_n^2}{\sqrt{(\lambda_1 - \lambda_n)^2 + 4\lambda_1^2 } + (\lambda_1 - \lambda_n)} \\
				= & \frac{2\lambda_n}{\sqrt{ (\kappa-1)^2 + 4\kappa^2 } + (\kappa - 1)} \\
				\geq & \frac{1}{\varphi} \frac{\lambda_n}{\kappa}.
			\end{aligned}
			\label{eq:lowerboundforC}
		\end{equation}
		Here, $\varphi \approx 1.6$ is the golden ratio and therefore $\delta \approx 0.6$. 
		
		\section*{Appendix: Proof of Proposition 4}
		
		First, note that in Algorithm 1, the spectrum of $\mathbf{S}^{(q)}$ is the same as that of $\mathbf{S}^{(0)} = \mathbf{S}$ because orthogonal transformations do not modify the eigenvalues (and therefore neither the condition number).
		
		We use the definition of $\mathbf{A}$ from Proposition 1. The extremal eigenvalues of the matrix $\mathbf{AS}^{(q)}\mathbf{A}$ are bounded by the extremal eigenvalues of $\mathbf{S}^{(q)}$ because $\mathbf{AS}^{(q)}\mathbf{A}$ is a principal sub-matrix of $\mathbf{S}^{(q)}$. Note that $\text{tr}(\mathbf{AS}^{(q)}) = \text{tr}(\mathbf{AS}^{(q)}\mathbf{A})$, because $\mathbf{A} = \mathbf{AA}$ followed by the cyclic permutation property of the trace $\text{tr}(\mathbf{AS}^{(q)})  = \text{tr}(\mathbf{AA} \mathbf{S}^{(q)})  = \text{tr}(\mathbf{AS}^{(q)}\mathbf{A})$, and $\mathbf{AS}^{(q)}\mathbf{A}$ is symmetric positive definite, it is, in fact, the principal sub-matrix of size $p \times p$ of $\mathbf{S}^{(q)}$. These results hold because of the particular choice $\mathbf{\alpha} = \mathbf{1}_{p \times 1}$.
		
		
		At each step of the proposed algorithm, the trace increases by the value of the scores we compute. For a single step of the algorithm, for the trace lower bound we have that:
		{\allowdisplaybreaks
			\begin{align*}
				\text{tr}(\mathbf{AS}^{(q)}) & =  \text{tr}(\mathbf{AG}_{i_q j_q}^T \mathbf{S}^{(q-1)} \mathbf{G}_{i_q j_q}) \\
				& \overset{Prop. 2}{=}   \text{tr}(\mathbf{AS}^{(q-1)}) + \mathcal{C}_{i_q j_q}   \\
				& \overset{Prop. 3}{\geq} \text{tr}(\mathbf{AS}^{(q-1)}) + \frac{\delta\lambda_n}{\kappa} \\
				& \overset{\eqref{eq:oneoverkappa}}{\geq} \text{tr}(\mathbf{AS}^{(q-1)}) + \frac{\delta \lambda_n^2}{\text{tr}(\mathbf{AS}^{(q-1)}) + (n-p)\lambda_1 - (n-1)\lambda_n} \\
				& \overset{\eqref{eq:eta}}{\geq}   \left(  1 + \frac{\delta \lambda_n^2}{p\lambda_1( n\lambda_1 - (n-1)\lambda_n ) } \right) \text{tr}(\mathbf{AS}^{(q-1)})  \\
				& \geq   \left(  1 + \frac{\delta \lambda_n^2}{pn\lambda_1^2 } \right) \text{tr}(\mathbf{AS}^{(q-1)})  \\
				& =  \left(  1 + \frac{\delta }{p n \kappa^2 } \right) \text{tr}(\mathbf{AS}^{(q-1)}).
			\end{align*}
		}
		We note that $\text{tr}(\mathbf{S}^{(q-1)}) = \text{tr}(\mathbf{AS}^{(q-1)}) +  \text{tr}((\mathbf{I}_n - \mathbf{A})\mathbf{S}^{(q-1)}) = \text{tr}(\mathbf{AS}^{(q-1)} \mathbf{A}) +  \text{tr}((\mathbf{I}_n - \mathbf{A})\mathbf{S}^{(q-1)}(\mathbf{I}_n - \mathbf{A}))$, then we used that
		\begin{equation}
			\begin{aligned}
				\frac{1}{\kappa}  = \frac{\lambda_n}{\lambda_1} & \overset{\eqref{eq:firsttraceinequality}}{\geq} \frac{\lambda_n}{\text{tr}(\mathbf{S}^{(q-1)})  - (n-1)\lambda_n} \\
				& = \frac{\lambda_n}{\text{tr}(\mathbf{AS}^{(q-1)}) +  \text{tr}((\mathbf{I}_n - \mathbf{A})\mathbf{S}^{(q-1)}) - (n-1)\lambda_n} \\
				& \overset{\eqref{eq:thirdtraceinequality}}{\geq} \frac{\lambda_n}{\text{tr}(\mathbf{AS}^{(q-1)}) + (n-p)\lambda_1 - (n-1)\lambda_n}.
			\end{aligned}
			\label{eq:oneoverkappa}
		\end{equation}
		Above, we have also upper-bounded the largest eigenvalue and we use the bounds for the trace quantities
		\begin{equation}
			\begin{aligned}
				(n-1)\lambda_n + \lambda_1 \leq &  \text{tr}(\mathbf{S}^{(q-1)}) \leq \lambda_n + (n-1)\lambda_1, \\
			\end{aligned}
			\label{eq:firsttraceinequality}
		\end{equation}
		\begin{equation}
			\begin{aligned}
				p\lambda_n\leq \sum_{t=1}^p \lambda_{n-t+1}  \leq  & \text{tr}(\mathbf{AS}^{(q-1)}\mathbf{A}) \leq \sum_{t=1}^p \lambda_t \leq  p\lambda_1,
			\end{aligned}
			\label{eq:secondtraceinequality}
		\end{equation}
		\begin{equation}
			\text{tr}((\mathbf{I}_n - \mathbf{A})\mathbf{S}^{(q-1)}) =  \text{tr}((\mathbf{I}_n - \mathbf{A})\mathbf{S}^{(q-1)}(\mathbf{I}_n - \mathbf{A})) \leq (n-p) \lambda_1.
			\label{eq:thirdtraceinequality}
		\end{equation}
		We also lower-bounded the function
		\begin{equation}
			x +  \frac{\delta \lambda_n^2 }{x + (n-p) \lambda_1 - (n-1)\lambda_n} \geq  (1 + \eta)x,
		\end{equation}
		for some constant $\eta$. To achieve the upper bound, we evaluate the function at the point $x = p\lambda_1 $, from the upper trace bound in \eqref{eq:secondtraceinequality}, because the fraction is a decreasing function of $x$ and therefore the minimum is given by the constant:
		\begin{equation}
			\eta  = \frac{\delta \lambda_n^2}{p\lambda_1(n\lambda_1 -(n-1)\lambda_n) }.
			\label{eq:eta}
		\end{equation}

		\section*{Appendix: Proof of Proposition 5}
		
		Let us define the following: $\mathcal{E}_q$ is the number of non-zero off-diagonal elements at step $q$ and the set of neighbors of an index $\mathcal{N}_q(i) = \{ j\ | \ (\mathbf{S}^{(q)})_{ij} \neq 0 \}$ at step $q$. We also define the symmetric difference operator between two sets, denoted as $\Delta$. Then, when at step $q$ we have chosen a transformation with indices $(i,j)$ we have that:
		\begin{equation}
			| \mathcal{E}_{q+1} | = | \mathcal{E}_{q} | + | \mathcal{N}_q(i) \Delta \mathcal{N}_q(j)  | + \mathbb{P}[ (\mathbf{S}^{(q)})_{ij} \neq 0 ].
		\end{equation}
		
		\noindent Based on this quantity, we can define the fill-in density to be:
		\begin{equation}
			p_q = \frac{ |\mathcal{E}_q | }{ N }, N = \frac{n(n-1)}{2}.
		\end{equation}
		
		The quantity $p_q$ can be interpreted as the probability that an off-diagonal element is different from zero.
		
		We have made the assumption that no exact cancellations take place in the calculations. To compute the size of the symmetric difference operator we consider fixed the indices $(i,j)$ and for a new index $k \notin \{ i, j\}$ a new element is added to the matrix when exactly one of the elements $(i,k)$ or $(j, k)$ is present. To quantify this situation, assuming independence, we have that:
		\begin{equation}
			\mathbb{P} [ \text{``exactly one''} ] \approx p_q(1-p_q) + (1-p_q)p_q = 2p_q(1-p_q),
		\end{equation}
		essentially assuming that $r_q \approx p_q^2$. Given this calculation, we have that for the symmetric difference that
		\begin{equation}
			\mathbb{E}[ \mathcal{N}_q(i) \Delta \mathcal{N}_q(j) \ | \ p_q ] \approx 2 (n-2) p_q(1-p_q),
		\end{equation}
		which in turn gives us the average density of new off-diagonal elements to be:
		\begin{equation}
			\mathbb{E}[  |\mathcal{E}_{q+1} | - | \mathcal{E}_q | \ | \ p_q  ] \approx \frac{4(n-2)}{n(n-1)} p_q (1-p_q) + \frac{2}{n(n-1)} (1-p_q).
		\end{equation}
		
		Keeping only the first dominant term of $O(1/n)$ we reach the final result:
		\begin{equation}
			p_q \approx \frac{1}{ 1 + C\exp(- q / O(n) )}.
		\end{equation}
		
		The independence assumption simplifies the result without affecting the qualitative conclusion, i.e., whether it is present or not, the sigmoid structure of the result remains unchanged. The exact result would be:
		\begin{equation}
			\begin{aligned}
				\mathbb{P} [ \text{``exactly one''} ] = & \mathbb{P}[ (\mathbf{S}^{(q)})_{ik} \neq 0 ] + \mathbb{P}[ (\mathbf{S}^{(q)})_{jk} \neq 0 ] \\
				& \ - 2\mathbb{P}[(\mathbf{S}^{(q)})_{ik} \neq 0 \text{ and } (\mathbf{S}^{(q)})_{jk} \neq 0 ] \\
				= & 2(p_q - r_q),
			\end{aligned}
		\end{equation}
		where we denoted by $r_q$ the expectation that two random rows share a non-zero entry at the same position. Although this quantity is trivial to compute in a concrete example, we were unable to develop a closed-form solution for the update of $r_q$.
		
		\section*{Acknowledgments}
		
		This research is supported by the project ``Romanian Hub for Artificial Intelligence - HRIA'', Smart Growth, Digitization and Financial Instruments Program, 2021-2027, MySMIS no. 334906.
		
		\bibliographystyle{elsarticle-num} 
		\bibliography{bib.bib, refs.bib}

		\vfill
		
	\end{document}